\documentclass[a4paper,10pt]{article}

\usepackage{amsmath,amssymb,amsfonts,amsthm}
\usepackage{mathptmx}
\allowdisplaybreaks[4]
\newcommand{\fz}{\frac}
\newcommand{\prz}[2]{ \frac{\partial{#1}}{\partial{#2}} }
\newcommand{\pz}{\partial}
\newcommand{\lA}{\langle}
\newcommand{\rA}{\rangle}
\newcommand{\rarrow}{\rightarrow}

\newcommand{\ol}{\overline}
\renewcommand{\epsilon}{\varepsilon}
\newcommand{\barO}{\bar{\Omega}}

\renewcommand{\Omega}{\varOmega}
\renewcommand{\Gamma}{\varGamma}
\renewcommand{\Psi}{\varPsi}
\renewcommand{\Pi}{\varPi}
\newtheorem{Def}{Definition}
\newtheorem{Thm}{Theorem}
\newtheorem{Prop}{Proposition}
\newtheorem{Lem}{Lemma}
\newtheorem{Hyp}{Hypothesis}
\newtheorem{Rmk}{Remark}

\providecommand{\keywords}[1]{\textit{Keywords:} #1}

\title{
Error estimates of stable and stabilized Lagrange-Galerkin schemes\\
for natural convection problems
}
\author{Hirofumi Notsu$^{1,}$\footnote{Corresponding author. E-mail: {\tt h.notsu@aoni.waseda.jp}} \ \ \  and \ Masahisa Tabata$^2$
\bigskip\\
$^1$ Waseda Institute for Advanced Study, Waseda University, Tokyo 169-8555, Japan \\
$^2$ Department of Mathematics, Waseda University, Tokyo 169-8555, Japan
}

\date{
\vspace*{-1em}
}
\begin{document}
\maketitle
\begin{abstract}
Optimal error estimates of stable and stabilized Lagrange-Galerkin (LG) schemes for natural convection problems are proved under a mild condition on time increment and mesh size.
The schemes maintain the common advantages of the LG method, i.e., robustness for convection-dominated problems and symmetry of the coefficient matrix of the system of linear equations.
We simply consider typical two sets of finite elements for the velocity, pressure and temperature, P2/P1/P2 and P1/P1/P1, which are employed by the stable and stabilized LG schemes, respectively.
The stabilized LG scheme has an additional advantage, a small number of degrees of freedom especially for three-dimensional problems.
The proof of the optimal error estimates is done by extending the arguments of the proofs of error estimates of stable and stabilized LG schemes for the Navier-Stokes equations in previous literature.
\smallskip\\
%
\end{abstract}
%
%
%
%
%
%
%
%
\section{Introduction}
In this paper we prove optimal error estimates of stable and stabilized Lagrange-Galerkin (LG) schemes for natural convection problems under a mild condition on the time increment and the mesh size.
\par
The LG method is the finite element method combined with an idea of the method of characteristics and has two common advantages, i.e., robustness for convection-dominated problems and symmetry of the coefficient matrix of the resultant system of linear equations.
There are several papers on error estimates of LG schemes for the Navier-Stokes equations, e.g., Pironneau~\cite{P-1982}, S\"uli~\cite{Suli-1988}, Boukir et al.~\cite{BMMR-1997} and Achdou and Guermond~\cite{AG-2000}, where the conventional inf-sup condition~\cite{GR-1986} is assumed to be satisfied and conforming stable elements, e.g., P2/P1 (Hood-Taylor) element~\cite{GR-1986}, are employed for them.
On the other hand, optimal error estimates of a stabilized LG scheme with one of the cheapest finite elements, P1/P1, have been proved in~\cite{NT-2015-M2AN}.
To the best of our knowledge, however, there is no proof of the optimal error estimates of an LG scheme for natural convection problems.
\par
In this paper we study stable and stabilized LG schemes of first-order in time for natural convection problems, and prove the optimal error estimates of the schemes.
We simply consider two typical sets of finite elements for the velocity, pressure and temperature, i.e., P2/P1/P2 and P1/P1/P1 elements, which are employed by the stable and stabilized LG schemes, respectively.
The schemes maintain the two common advantages of the LG method mentioned above.
We note that the pair of the velocity and pressure and the temperature are solved alternatively in the schemes, and that both of the resulting matrices are symmetric.
The stabilized LG scheme has an additional advantage, a small number of degrees of freedom, which leads to efficient computation especially in three-dimensions.
\par
The proof of the optimal error estimates is performed by extending the arguments of the proofs of error estimates of stable and stabilized LG schemes for the Navier-Stokes equations established in~\cite{NT-2015-M2AN,Suli-1988}.
The essential part of the argument of the proof for the stable LG scheme~\cite{Suli-1988} is as follows.
If the equation in the error~$(e_{uh}, e_{ph})$ of the form
\begin{align}
\int_\Omega \fz{1}{\Delta t} (e_{uh}^n-e_{uh}^{n-1}) \cdot v_h~dx + \int_\Omega D(e_{uh}^n) : D(v_h)~dx - \int_\Omega (\nabla\cdot v_h) e_{ph}^n~dx - \int_\Omega (\nabla\cdot e_{uh}^n) q_h~dx = \int_\Omega R_{uh}^n \cdot v_h~dx, \notag\\
\forall (v_h, q_h)\in V_h\times Q_h,
\label{eq:error_NS_introduction}
\end{align}
with the estimate of the remainder~$R_{uh}^n$ in the $L^2(\Omega)^d$-norm
\begin{align*}
\|R_{uh}^n\|_{L^2(\Omega)^d} \le c (\|u_h^{n-1}\|_{L^\infty(\Omega)}+1) (\|e_{uh}^{n-1}\|_{H^1(\Omega)^d} + \Delta t + h^k)
\end{align*}
is obtained, we can show the conditional stability and the optimal error estimates by mathematical induction,
where 
$\Delta t$ is a time increment, 
$t^n = n\Delta t$ is a time at step~$n$, 
$h$ is a mesh size,
$d$ is the space dimension,
$k$ is a positive integer depending on the choice of finite element spaces~$V_h$ for the velocity and $Q_h$ for the pressure,
$D(v)$ is the strain-rate tensor with respect to the velocity~$v$,
$(e_{uh}, e_{ph})=\{(e_{uh}^n, e_{ph}^n)\}_n$ is a set of functions of error for the velocity and pressure defined by $(e_{uh}^n, e_{ph}^n) = (u_h^n, p_h^n) - (\hat{u}_h, \hat{p}_h)(t^n)$, 
$(u_h, p_h)=\{(u_h^n, p_h^n)\}_n$ is the solution of the scheme,
and $(\hat{u}_h, \hat{p}_h)=\{(\hat{u}_h, \hat{p}_h)(t)\}_t$ is a Stokes projection of the exact solution~$(u,p) = \{(u, p)(t)\}_t$.
The key issue of the argument is that the value~$\|u_h^{n-1}\|_{L^\infty(\Omega)^d}$, $L^\infty(\Omega)^d$-norm of the numerical velocity at the previous time step, can be employed for the estimate of~$\|R_{uh}^n\|_{L^2(\Omega)^d}$.
\par
The natural convection problem consists of a system of equations of velocity, pressure and temperature.
For the error estimates of the stable LG scheme for the problem, we can extend the argument above to the scheme by considering a set of functions of error~$(e_{uh}, e_{ph}, e_{\theta h})=\{(e_{uh}^n, e_{ph}^n, e_{\theta h}^n)\}_n$ defined by $(e_{uh}^n, e_{ph}^n, e_{\theta h}^n)=(u_h^n, p_h^n, \theta_h^n) - (\hat{u}_h, \hat{p}_h, \hat{\theta}_h)(t^n)$ for the solution of the scheme~$(u_h, p_h, \theta_h)=\{(u_h^n, p_h^n, \theta_h^n)\}_n$ and a Stokes-Poisson projection~$(\hat{u}_h, \hat{p}_h, \hat{\theta}_h) = \{(\hat{u}_h, \hat{p}_h, \hat{\theta}_h)(t)\}_t$ of the exact solution~$(u, p, \theta)=\{(u, p, \theta)(t)\}_t$.
In fact, we derive a corresponding equation in error~$(e_{uh}, e_{ph}, e_{\theta h})$ to~\eqref{eq:error_NS_introduction}
\begin{align*}
\int_\Omega \fz{1}{\Delta t} \Bigl\{ (e_{uh}^n-e_{uh}^{n-1}) \cdot v_h + (e_{\theta h}^n-e_{\theta h}^{n-1}) \psi_h \Bigr\}~dx 
+ \int_\Omega \Bigl\{ D(e_{uh}^n) : D(v_h) + \nabla e_{\theta h}^n \cdot \nabla \psi_h \Bigr\}~dx 
- \int_\Omega (\nabla\cdot v_h) e_{ph}^n~dx \qquad\\
- \int_\Omega (\nabla\cdot e_{uh}^n) q_h~dx 
= \int_\Omega (R_{uh}^n \cdot v_h + R_{\theta h}^n \psi_h)~dx,\quad \forall (v_h, q_h, \psi_h) \in V_h\times Q_h\times \Psi_h,
\end{align*}
with the estimates of the remainders~$R_{uh}^n$ and $R_{\theta h}^n$
\begin{align*}
\|R_{uh}^n\|_{L^2(\Omega)^d},\ \|R_{\theta h}^n\|_{L^2(\Omega)} \le c (\|u_h^{n-1}\|_{L^\infty(\Omega)^d}+1) (\|e_{uh}^{n-1}\|_{H^1(\Omega)^d}+\|e_{\theta h}^{n-1}\|_{H^1(\Omega)} + \Delta t + h^k),
\end{align*}
where $V_h$, $Q_h$ and $\Psi_h$ are finite element spaces for the velocity, pressure and temperature, respectively.
Then, we prove the optimal error estimates by mathematical induction.
From the discussions of stabilized LG schemes for the Oseen and the Navier-Stokes equations in~\cite{NT-2015-JSC,NT-2015-M2AN} it is not difficult to extend the argument to a corresponding stabilized LG scheme for natural convection problems.
Thus, we prove the optimal error estimates of both stable and stabilized LG schemes in this paper.
\par
This paper is organized as follows.
Stable and stabilized LG schemes for natural convection problems are presented in Section~\ref{sec:schemes}.
The main results, the conditional stability and the optimal error estimates, are stated in Section~\ref{sec:main_results}, and they are proved in Section~\ref{sec:proofs}.
The conclusions are given in Section~\ref{sec:conclusions}.
%
%
\section{Lagrange-Galerkin schemes}\label{sec:schemes}
The function spaces and the notation to be used throughout the paper are as follows.
Let $\Omega$ be a bounded domain in $\mathbb{R}^d (d=2,3)$, $\Gamma\equiv\pz\Omega$ the boundary of $\Omega$, and $T$ a positive constant.
For $m \in \mathbb{N}\cup \{0\}$ and $p\in [1,\infty]$ we use the {\rm Sobolev} spaces $W^{m,p}(\Omega)$, $W^{1,\infty}_0(\Omega)$, $H^m(\Omega) \, (=W^{m,2}(\Omega))$ and $H^1_0(\Omega)$.
For any normed space $S$ with norm $\|\cdot\|_S$, we define function spaces $C([0,T]; S)$ and $H^m(0,T; S)$ consisting of $S$-valued functions in $C([0,T])$ and $H^m(0,T)$, respectively.
We use the same notation $(\cdot, \cdot)$ to represent the $L^2(\Omega)$ inner product for scalar-, vector- and matrix-valued functions.
The dual pairing between $S$ and the dual space $S^\prime$ is denoted by $\lA\cdot, \cdot\rA$.
The norms on $W^{m,p}(\Omega)^d$ and $H^m(\Omega)^d$ are simply denoted by
$\|\cdot\|_{m,p}$ and $\|\cdot\|_m \, (= \|\cdot\|_{m,2})$, respectively,
and the notation~$\|\cdot\|_m$ is employed not only for vector-valued functions but also for scalar-valued ones.
We often omit $[0,T]$, $\Omega$ and/or $d$ if there is no confusion, e.g., we shall write $C(L^\infty)$ in place of $C([0,T]; L^\infty(\Omega)^d)$.
For $t_0$ and $t_1\in\mathbb{R}$ we introduce the function space,
\begin{align*}
Z^m(t_0, t_1) \equiv \Bigl\{ \psi \in H^j(t_0, t_1; H^{m-j}(\Omega));~j=0,\cdots,m,\ \|\psi\|_{Z^m(t_0, t_1)} < \infty \Bigr\}
\ \ \mbox{with}\ \ 
\|\psi\|_{Z^m(t_0, t_1)} \equiv \biggl\{ \sum_{j=0}^m \|\psi\|_{H^j(t_0,t_1; H^{m-j}(\Omega))}^2 \biggr\}^{1/2},
\end{align*}
and set $Z^m \equiv Z^m(0, T)$.
Hereafter we use special notations, $X\equiv H^1(\Omega)$, $Y\equiv X^d$, $M\equiv L^2(\Omega)$, $V\equiv H^1_0(\Omega)^d$, $Q \equiv L^2_0(\Omega) \equiv \{ q \in M;~(q, 1)=0\}$, $\Psi \equiv H^1_0(\Omega)$ and $\mathbb{H}^{m+1} \equiv H^{m+1}(\Omega)^d \times H^m(\Omega) \times H^{m+1}(\Omega)$, where $\|\cdot\|_V \equiv \|\cdot\|_1$, $\|\cdot\|_Q \equiv \|\cdot\|_0$ and $\|\cdot\|_\Psi \equiv \|\cdot\|_1$.
\par
We consider the following natural convection problem; find $(u, p, \theta):$ $\Omega\times (0, T) \rarrow \mathbb{R}^d\times\mathbb{R}\times\mathbb{R}$ such that
\begin{subequations}\label{prob:NCP}
\begin{align}
\frac{Du}{Dt} - \nabla \cdot [ 2\nu D(u) ] + \nabla p -\theta \beta & = f_u & & \mbox{in}\ \ \Omega\times (0,T), \label{eq:NCP_1}\\
\nabla \cdot u & = 0 & & \mbox{in}\ \ \Omega\times (0,T), \\
\frac{D\theta}{Dt} - \kappa\Delta \theta & = f_\theta & &\mbox{in}\ \ \Omega\times (0,T), \\
u = 0, \quad \theta & = 0 && \mbox{on}\ \ \Gamma\times (0,T),\\
u = u^0, \quad \theta & = \theta^0 & & \mbox{in}\ \ \Omega,\ \mbox{at}\ t=0,
\end{align}
\end{subequations}
where $u$ is the velocity, $p$ is the pressure, $\theta$ is the temperature, $f_u \in C([0,T]; L^2(\Omega)^d)$ is an external force, $f_\theta \in C([0,T]; L^2(\Omega))$ is a heat source, $\nu>0$ is a viscosity, $\kappa>0$ is a thermal conductivity, $\beta\in C([0,T]; L^\infty(\Omega)^d)$ is a generalized thermal expansion, $u^0 \in V$ is an initial velocity, $\theta^0 \in \Psi$ is an initial temperature, $D(u) \equiv (1/2)[\nabla u + (\nabla u)^T]$ is the strain-rate tensor and $D/Dt$ is the material derivative defined by $D/Dt \equiv \pz/\pz t + u\cdot\nabla$.
\begin{Rmk}
The last term of the left-hand side of~\eqref{eq:NCP_1}, $-\theta \beta$, is a generalized expression of the {\rm Boussinesq} approximation.
Indeed, buoyancy-driven flows can be treated by setting $\beta (x,t) = (0, 0, \beta_3(x,t))$, cf.~\cite{LB-1999,TT-2005}.
\end{Rmk}
\par
We have the weak formulation of~\eqref{prob:NCP}; find $(u, p, \theta): (0, T)\to V\times Q\times \Psi$ such that, for $t\in (0,T)$,
\begin{align}
\biggl( \fz{Du}{Dt}(t), v \biggr)+\biggl( \fz{D\theta}{Dt}(t), \psi \biggr) + \mathcal{A} \bigl( (u,p)(t), (v,q) \bigr) +  a_\theta (\theta(t),\psi) - \bigl(\theta(t)\beta(t), v\bigr) = \bigl( f_u (t), v \bigr) + \bigl( f_\theta (t), \psi \bigr), \notag\\
\forall (v, q, \psi)\in V\times Q\times \Psi,
\label{eq:NCP_weak}
\end{align}
with $u(0) = u^0$ and $\theta(0) = \theta^0$, where $\mathcal{A}$ and $a_\theta$ are bilinear forms on $(V\times Q) \times (V\times Q)$ and $\Psi \times \Psi$ defined by
\begin{align*}
&\mathcal{A} \bigl( (u,p), (v,q) \bigr) \equiv a_u (u,v)+b(v,p)+b(u,q),\\
& a_u (u,v) \equiv 2 \nu \bigl( D(u), D(v) \bigr),
\quad
b(v,q) \equiv -( \nabla\cdot v, q ),
\quad
a_\theta (\theta,\psi) \equiv \kappa ( \nabla\theta, \nabla\psi ).
\end{align*}
\par
We introduce the method of characteristics.
Let $\Delta t$ be a time increment, $N_T\equiv \lfloor T/\Delta t \rfloor$ a total number of time steps and $t^n \equiv n\Delta t$ for $n=0,\cdots,N_T$.
Let $\phi$ and $\psi$ be functions defined in $\Omega\times (0,T)$ and $\Omega$, respectively.
We generally denote $\phi(\cdot,t^n)$ by $\phi^n$.
Let $X=X(\cdot; x, t^n): (0, T)\to\mathbb{R}^d$ be the solution of the system of ordinary differential equations,
\begin{align*}
\frac{dX}{dt}(t) = u\bigl( X(t), t \bigr),
\end{align*}
with an initial condition~$X(t^n)=x$.
Then, we obtain a first-order approximation of $D\phi/Dt$ at $(x,t^n)$ as follows:
\begin{align*}
\fz{D\phi}{Dt} (x, t^n) & = \fz{D\phi}{Dt} (X(t), t)\Bigr|_{t=t^n} = \frac{d}{dt} \phi \bigl( X(t), t \bigr)\Bigr|_{t=t^n}
= \frac{\phi^n - \phi^{n-1}\circ X_1(u^{n-1},\Delta t)}{\Delta t} (x) + O(\Delta t),
\end{align*}
where the symbol $\circ$ stands for the composition of functions $(\psi \circ v) (x) \equiv  \psi\bigl( v(x) \bigr)$ for $v: \Omega \to \Omega$, and $X_1(w,\Delta t): \Omega\to\mathbb{R}^d$ is a mapping defined by
\begin{align*}
X_1(w,\Delta t)(x) \equiv x - w(x)\Delta t
\end{align*}
for $w: \Omega \to \mathbb{R}^d$.
The next proposition presents sufficient conditions to ensure that all \textit{upwind} points by $X_1(w,\Delta t)$ are in $\Omega$ and that its Jacobian~$J=J(w,\Delta t) \equiv \det ( \pz X_1(w,\Delta t) / \pz x )$ is around $1$.
\begin{Prop}[{\cite{RT-2002,TU-2015_NS}}]\label{prop:RT_TU}
(i)~Let $w\in W^{1,\infty}_0(\Omega)^d$ be a given velocity.
Then, under the condition $\Delta t |w|_{1,\infty} < 1$, $X_1(w,\Delta t): \Omega \to \Omega$ is bijective.
\ \ 
(ii)~Furthermore, under the condition $\Delta t |w|_{1,\infty} \le 1/4$, the estimate $1/2 \le J(w,\Delta t) \le 3/2$ holds.
\end{Prop}
\par
For the sake of simplicity we assume that $\Omega$ is a polygonal $(d=2)$ or polyhedral $(d=3)$ domain.
Let $\mathcal{T}_h=\{K\}$ be a triangulation of $\barO\ (= \bigcup_{K\in\mathcal{T}_h} K )$, $h_K$ a diameter of $K\in\mathcal{T}_h$, and $h\equiv\max_{K\in\mathcal{T}_h}h_K$ the maximum element size.
Throughout this paper we consider a regular family of triangulations $\{\mathcal{T}_h\}_{h\downarrow 0}$ with the inverse assumption~\cite{C-1978}, i.e., there exists a positive constant $\alpha_0$ independent of $h$ such that $h/h_K \le \alpha_0$ for any $K\in \mathcal{T}_h$ and $h$.
We define two sets of finite element spaces, $(V_h, Q_h, \Psi_h)$, depending on $k$ as follows.
Let
\begin{align}
k = 1\mbox{\ \ or\ \ } 2
\label{def:k}
\end{align}
be a fixed integer.
Let $X_h$, $M_h$ and $Y_h$ be finite element spaces defined by
\begin{align*}
X_h \equiv \{ v_h\in C(\barO);\ v_{h|K} \in P_k (K),\ \forall K \in\mathcal{T}_h \} \subset X,
\quad
Y_h \equiv X_h^d \subset Y,
\quad
M_h \equiv \{ q_h\in C(\barO);\ q_{h|K} \in P_1 (K), \ \forall K \in\mathcal{T}_h \} \subset M,
\end{align*}
where for $i=1, 2$, $P_i (K)$ is the space of polynomial functions of degree~$i$ on $K\in\mathcal{T}_h$.
We set $V_h \equiv V \cap Y_h$, $Q_h \equiv Q \cap M_h$ and $\Psi_h \equiv \Psi \cap X_h$.
\par
Let $(\cdot,\cdot)_K$ be the $L^2(K)^d$ inner product, and $\mathcal{C}_h$ and $\mathcal{A}_h$ bilinear forms on $X \times X$ and $(Y \times X) \times (Y \times X)$ defined by
\begin{align*}
\mathcal{C}_h (p, q) \equiv \sum_{K\in\mathcal{T}_h} h_K^2 ( \nabla p, \nabla q )_K,
\qquad
\mathcal{A}_h \bigl( (u,p), (v,q) \bigr) \equiv \mathcal{A} \bigl( (u,p), (v,q) \bigr) - \delta_0\mathcal{C}_h(p, q),
\end{align*}
respectively, where $\delta_0$ is defined by
\begin{align*}
\delta_0 & \equiv 
\left\{
\begin{aligned}
& 0 && (k = 2)\\
& 1 && (k = 1)
\end{aligned}
\right..
\end{align*}
\begin{Rmk}
The well-known stable element P2/P1~($k = 2$) satisfies the convectional inf-sup condition,
\begin{align}
\inf_{q_h \in Q_h}\sup_{v_h\in V_h} \fz{b(v_h, q_h)}{\|v_h\|_{V} \|q_h\|_{Q}} \ge \beta^\ast,
\label{ieq:conventional_inf_sup}
\end{align}
and the cheap equal-order finite element P1/P1~($k = 1$) satisfies a general version of~\eqref{ieq:conventional_inf_sup}, cf.~\cite{FS-1991},
\begin{align}
\inf_{(u_h, p_h) \in V_h\times Q_h}\sup_{(v_h, q_h) \in V_h\times Q_h} \fz{\mathcal{A}_h\bigl( (u_h, p_h), (v_h, q_h)\bigr)}{\|(u_h,p_h)\|_{V\times Q} \|(v_h, q_h)\|_{V\times Q}} \ge \gamma^\ast,
\label{ieq:generalized_inf_sup}
\end{align}
where $\beta^\ast$ and $\gamma^\ast$ are positive constants independent of $h$.
Since \eqref{ieq:conventional_inf_sup} implies \eqref{ieq:generalized_inf_sup}, \eqref{ieq:generalized_inf_sup} is satisfied for both cases~$k=1$ and~$2$.
\end{Rmk}
\begin{Rmk}
We simply set $\delta_0=1$ in the case~$k=1$.
For more discussion on the choice of $\delta_0$ see, e.g.,~\cite{Tez-2007,TMRS-1992}.
\end{Rmk}
\par
Suppose that the pair of approximate initial values, $(u_h^0, \theta_h^0) \in V_h \times \Psi_h$, of $(u^0, \theta^0) \in V \times \Psi$ is given.
The LG scheme for~\eqref{prob:NCP} is to find $\{(u_h^n, p_h^n, \theta_h^n)\}_{n=1}^{N_T}\subset V_h \times Q_h\times \Psi_h$ such that, for $n=1,\cdots, N_T$,
\begin{align}
\biggl(\fz{u_h^n-u_h^{n-1}\circ X_1(u_h^{n-1},\Delta t)}{\Delta t}, v_h\biggr) 
+\biggl(\fz{\theta_h^n-\theta_h^{n-1}\circ X_1(u_h^{n-1},\Delta t)}{\Delta t}, \psi_h\biggr)
+ \mathcal{A}_h \bigl( (u_h^n, p_h^n), (v_h, q_h) \bigr) + a_\theta (\theta_h^n, \psi_h) - (\theta_h^{n-1}\beta^n, v_h) \notag\\
= (f_u^n, v_h)+(f_\theta^n, \psi_h),\quad \forall (v_h, q_h, \psi_h) \in V_h\times Q_h\times \Psi_h.
\label{scheme_NCP}
\end{align}
Scheme~\eqref{scheme_NCP} is equivalent to
\begin{subequations}\label{scheme_NCP_equivalent}
\begin{align}
\biggl(\fz{u_h^n-u_h^{n-1}\circ X_1(u_h^{n-1},\Delta t)}{\Delta t}, v_h\biggr) 
+ a_u (u_h^n, v_h)+b(v_h, p_h^n) -(\theta_h^{n-1}\beta^n, v_h) & = (f_u^n, v_h), & \forall v_h & \in V_h,\label{scheme_NCP_equivalent:eq1}\\
b(u_h^n, q_h) - \delta_0\mathcal{C}_h(p_h^n, q_h) & = 0, & \forall q_h & \in Q_h,\label{scheme_NCP_equivalent:eq2}\\
\biggl(\fz{\theta_h^n-\theta_h^{n-1}\circ X_1(u_h^{n-1},\Delta t)}{\Delta t}, \psi_h\biggr)
+ a_\theta (\theta_h^n, \psi_h) & = (f_\theta^n, \psi_h), & \forall \psi_h & \in \Psi_h.\label{scheme_NCP_equivalent:eq3}
\end{align}
\end{subequations}
\par
In the following we often call scheme~\eqref{scheme_NCP} with P2/P1/P2~($k=2$) or P1/P1/P1~($k=1$) the stable or stabilized scheme, respectively.
\begin{Rmk}\label{rmk:scheme}
Suppose that a pair~$(u_h^{n-1}, \theta_h^{n-1}) \in V_h \times \Psi_h$ is given.
Under the condition~$\Delta t |u_h^{n-1}|_{1,\infty} < 1$ we have $X_1(u_h^{n-1},\Delta t)(\Omega)$ $=\Omega$ by Proposition~\ref{prop:RT_TU}.
Since $\theta_h^{n-1}$ is employed for the last term of the left-hand side of~\eqref{scheme_NCP_equivalent:eq1}, we separately get $(u_h^n,p_h^n)\in V_h\times Q_h$ and $\theta_h^n\in \Psi_h$ from~\eqref{scheme_NCP_equivalent:eq1}-\eqref{scheme_NCP_equivalent:eq2} and~\eqref{scheme_NCP_equivalent:eq3}, respectively, where both resulting matrices are invertible and symmetric.
The invertibility of the matrix of~\eqref{scheme_NCP_equivalent:eq1}-\eqref{scheme_NCP_equivalent:eq2} is assured by virtue of the inf-sup condition~\eqref{ieq:conventional_inf_sup} and {\rm Brezzi-Pitk\"aranta}'s stabilization~$\mathcal{C}_h$ for the stable and stabilized schemes, respectively.
Thus, we get the unique solution $(u_h^n, p_h^n, \theta_h^n)\in V_h\times Q_h\times \Psi_h$ of~\eqref{scheme_NCP}.
\end{Rmk}
%
%
%
%
%
%
%
%
%
%
\section{Main results}\label{sec:main_results}
In this section we state the main results, conditional stability and optimal error estimates for scheme~\eqref{scheme_NCP}, which are proved in Section~\ref{sec:proofs}.
We use the following norms and a seminorm,
\begin{align*}
\|u\|_{\ell^{\infty}(S)} &\equiv \max_{n=0,\cdots, N_T} \|u^n\|_{S},
&
\|u\|_{\ell^{2}_m(S)} &\equiv \biggl\{ \Delta t\sum_{n=1}^{m} \|u^n\|_S^2 \biggr\}^{1/2},
&
\|u\|_{\ell^{2}(S)} &\equiv \|u\|_{\ell^{2}_{N_T}(S)},
&
|p|_h & \equiv \mathcal{C}_h (p, p)^{1/2},
\end{align*}
for $m\in\{1,\cdots,N_T\}$ and $S=L^\infty(\Omega)$, $L^2(\Omega)$ and $H^1(\Omega)$.
Let $\ol{D}_{\Delta t}$ be the backward difference operator defined by~$\ol{D}_{\Delta t} u^n \equiv (u^n - u^{n-1})/\Delta t$.
\begin{Hyp}\label{hyp:regularity}
The solution $(u, p, \theta)$ of~\eqref{eq:NCP_weak} satisfies
$u\in C([0,T]; W^{1,\infty}(\Omega)^d)\cap (Z^2)^d \cap H^1(0,T; V\cap H^{k+1} (\Omega)^d)$, $p\in H^1(0,T;Q\cap H^k (\Omega))$ and~$\theta \in C([0,T]; W^{1,\infty}(\Omega))\cap Z^2 \cap H^1(0,T; \Psi\cap H^{k+1} (\Omega))$ for the integer $k$ in~\eqref{def:k}.
\end{Hyp}
\begin{Def}[Stokes-Poisson projection]\label{def:StokesPoissonPrj}
For $(w,r,\phi) \in V \times Q \times \Psi$ we define the Stokes-Poisson projection $(\hat{w}_h, \hat{r}_h, \hat{\phi}_h) \in V_h\times Q_h \times \Psi_h$ of $(w, r, \phi)$ by
\begin{align}
\mathcal{A}_h\bigl( (\hat{w}_h, \hat{r}_h), (v_h, q_h) \bigr) + a_\theta (\hat{\phi}_h, \psi_h)= \mathcal{A} \bigl( (w, r), (v_h, q_h) \bigr) + a_\theta (\phi, \psi_h), \quad \forall (v_h, q_h,\psi_h) \in V_h\times Q_h \times \Psi_h.
\label{eq:StokesPoissonPrj}
\end{align}
\end{Def}
\par
Since the Stokes-Poisson projection is well-defined, there exists a projection operator $\Pi_h^{\rm SP}: V \times Q \times \Psi \to V_h\times Q_h \times \Psi_h$ defined by $\Pi_h^{\rm SP} (w,r,\phi) = (\hat{w}_h, \hat{r}_h, \hat{\phi}_h)$.
We denote the $i$-th component of $\Pi_h^{\rm SP} (w,r,\phi)$ by $[\Pi_h^{\rm SP} (w,r,\phi)]_i$ for $i=1,2,3$ and the pair of the first and third components~$(\hat{w}_h, \hat{\phi}_h)=([\Pi_h^{\rm SP} (w,r,\phi)]_1, [\Pi_h^{\rm SP} (w,r,\phi)]_3)$ by $[\Pi_h^{\rm SP} (w,r,\phi)]_{1,3}$ simply.
\par
We state the main results.
\begin{Thm}\label{thm:main_results}
Let $k$ be the integer in~\eqref{def:k}.
Suppose Hypothesis~\ref{hyp:regularity} holds.
Then, there exist positive constants $h_0$ and $c_0$ independent of $h$ and $\Delta t$ such that, for any pair~$(h, \Delta t)$,
\begin{align}
h\in (0,h_0],\quad \Delta t\le c_0 h^{d/4},
\label{condition:h_dt}
\end{align}
the following hold.
\smallskip\\
(i)~Scheme~\eqref{scheme_NCP} with $(u_h^0,\theta_h^0)=[\Pi_h^{\rm SP}(u^0, 0,\theta^0)]_{1,3}$ has a unique solution~$(u_h, p_h, \theta_h)=\{(u_h^n, p_h^n, \theta_h^n)\}_{n=1}^{N_T}\subset V_h\times Q_h\times \Psi_h$.
\medskip\\
(ii)~It holds that
\begin{align}
\|u_h\|_{\ell^\infty(L^\infty)}\le \|u\|_{C(L^\infty)}+1.
\label{ieq:stability}
\end{align}
(iii)~There exists a positive constant~$c_\dagger$ independent of $h$ and $\Delta t$ such that
\begin{align}
\|u_h-u\|_{\ell^\infty(H^1)},\ \ \Bigl\|\ol{D}_{\Delta t}u_h-\prz{u}{t}\Bigr\|_{\ell^2(L^2)},\ \ \|p_h-p\|_{\ell^2(L^2)},\ \ \|\theta_h-\theta\|_{\ell^\infty(H^1)},\ \ \Bigl\|\ol{D}_{\Delta t}\theta_h-\prz{\theta}{t}\Bigr\|_{\ell^2(L^2)} \le c_\dagger (\Delta t + h^k).
\label{ieq:main_results}
\end{align}
\end{Thm}
\begin{Hyp}\label{hyp:L2}
The Stokes and Poisson problems are both regular, i.e., for any $g_S \in L^2(\Omega)^d$ and $g_P \in L^2(\Omega)$, the solution $(w, r)\in V\times Q$ of the Stokes problem,
\begin{align*}
a_u(w,v)+b(v,r)+b(w,q) = (g_S, v), \quad \forall (v, q) \in V\times Q,
\end{align*}
and the solution $\phi \in \Psi$ of the Poisson problem,
\begin{align*}
a_\theta(\phi, \psi) = (g_P,\psi), \quad \forall \psi \in \Psi,
\end{align*}
belong to $H^2(\Omega)^d\times H^1(\Omega)$ and $H^2(\Omega)$, respectively,
and the estimates
\begin{align*}
\|w\|_2 + \|r\|_1 \le c_R \|g_S\|_0,\qquad \|\phi\|_2 \le \bar{c}_R \|g_P\|_0,
\end{align*}
hold, where $c_R$ and $\bar{c}_R$ are positive constants independent of $g_S$, $g_P$, $w$, $r$ and $\phi$.
\end{Hyp}
\begin{Thm}\label{thm:main_results_L2}
Let $k$ be the integer in~\eqref{def:k}.
Suppose Hypotheses~\ref{hyp:regularity} and~\ref{hyp:L2} hold.
Then, there exists a positive constant~$c_\ddagger$ independent of $h$ and $\Delta t$ such that
\begin{align}
\|u_h-u\|_{\ell^\infty(L^2)},\ \ \|\theta_h-\theta\|_{\ell^\infty(L^2)} \le c_\ddagger (\Delta t + h^{k+1}),
\label{ieq:L2}
\end{align}
where $u_h$ and $\theta_h$ are the first and third components of the solution of~\eqref{scheme_NCP} of Theorem~\ref{thm:main_results}-(i).
\end{Thm}
\begin{Rmk}
We prepare the approximate initial value by the Stokes-Poisson projection, i.e., $(u_h^0,\theta_h^0)=[\Pi_h^{\rm SP}(u^0, 0,\theta^0)]_{1,3}$, which does not lose any convergence orders in the two theorems above.
\end{Rmk}
\begin{Rmk}
Hypothesis~\ref{hyp:L2} holds, e.g., if $\Omega$ is convex in $\mathbb{R}^2$, cf.~\cite{GR-1986}.
\end{Rmk}
\begin{Rmk}\label{rmk:generalization}
(i)~To general stable {\rm Hood-Taylor}(-type) elements P$\ell$/P($\ell$-1)/P$\ell$, $\ell>2$, with~\eqref{ieq:conventional_inf_sup}, the results of Theorems~\ref{thm:main_results} and~\ref{thm:main_results_L2} can be extended, where the error estimates are of orders $O(\Delta t+h^\ell)$ and $O(\Delta t+h^{\ell+1})$, respectively.
\ \ 
(ii)~For the stable mini(-type) element P1+/P1/P1, cf.~\cite{GR-1986}, we can prove~\eqref{ieq:main_results} and~\eqref{ieq:L2} with $k=1$ for $\delta_0=0$.
\ \ 
(iii)~For generalization of the stabilized {\rm LG} scheme with, e.g., P$\ell$/P$\ell$/P$\ell$ element, $\ell>1$, a modification of the bilinear form~$\mathcal{C}$ is required in order to obtain optimal error estimates, cf.~\cite{Bur-2008}.
\end{Rmk}
%
%
%
%
%
%
%
%
%
\section{Proofs}\label{sec:proofs}
In this section Theorems~\ref{thm:main_results} and~\ref{thm:main_results_L2} are proved.
We use generic positive constants~$c$, $c_u$, $c_\theta$, $c_{(u,\theta)}$ and $c_{(u,p,\theta)}$ independent of $h$ and $\Delta t$.
$c_u$, $c_\theta$, $c_{(u,\theta)}$ and $c_{(u,p,\theta)}$ depend on $u$, $\theta$, $(u,\theta)$ and $(u,p,\theta)$, respectively.
\subsection{Preparations}\label{subsec:preparations}
We prepare lemmas and a proposition, which are directly used in our proofs.
\begin{Lem}[\cite{DL-1976}]\label{lem:Korn}
Let $\Omega$ be a bounded domain with a Lipschitz-continuous boundary.
Then, there exist positive constants $\alpha_1$ and $\bar{\alpha}_1$ and the following inequalities hold.
\begin{align*}
\|D(v)\|_0 & \le \|v\|_1 \le \alpha_1 \|D(v)\|_0, & \forall v & \in H^1_0(\Omega)^d,
&
\|\nabla \psi\|_0 & \le \|\psi\|_1 \le \bar{\alpha}_1 \|\nabla\psi\|_0, & \forall \psi & \in H^1_0(\Omega).
\end{align*}
\end{Lem}
\begin{Lem}[\cite{C-1978}]\label{lem:inverse_inequality}
There exist positive constants $\alpha_{2i}$, $i=0,\cdots,4$, independent of $h$ and the following hold.
\begin{align*}
|q_h|_h & \le \alpha_{20} \|q_h\|_0, & \forall q_h & \in Q_h,
&
\|v_h\|_{0,\infty} & \le \alpha_{21} h^{-d/6}\|v_h\|_1, & \forall v_h & \in V_h,
\\
\|v_h\|_{1,\infty} & \le \alpha_{22} h^{-d/2}\|v_h\|_1, & \forall v_h & \in V_h,
&
\|\Pi_h v\|_{0,\infty} & \le \|v\|_{0,\infty}, & \forall v & \in C(\bar{\Omega})^d,
\\
\|\Pi_h v\|_{1,\infty} & \le \alpha_{23} \|v\|_{1,\infty}, & \forall v & \in W^{1,\infty}(\Omega)^d,
&
\|\Pi_h v - v\|_1 & \le \alpha_{24} h^k \|v\|_{k+1}, & \forall v & \in H^{k+1}(\Omega)^d,
\end{align*}
where $\Pi_h: C(\bar{\Omega})^d \to Y_h$ is the {\rm Lagrange} interpolation operator.
\end{Lem}
\begin{Rmk}\label{rmk:interpolation_property}
We note $\alpha_{23} \ge 1$.
\end{Rmk}
\begin{Prop}\label{prop:StokesPoissonPrj}
(i)~Let $k$ be the integer in~\eqref{def:k}.
Suppose $(w,r,\phi) \in (V\times Q \times \Psi) \cap \mathbb{H}^{k +1}$.
Then, there exist positive constants~$\alpha_{31}$ and~$\bar{\alpha}_{31}$ independent of $h$ such that the Stokes-Poisson projection~$(\hat{w}_h, \hat{r}_h, \hat{\phi}_h) \in V_h\times Q_h \times \Psi_h$ of $(w, r, \phi)$ by~\eqref{eq:StokesPoissonPrj} satisfies the estimates,
\begin{subequations}
\begin{align}
\|\hat{w}_h-w\|_1,\ \ \|\hat{r}_h-r\|_0,\ \ \delta_0|\hat{r}_h-r|_h \le \alpha_{31} h^{k} \|(w,r)\|_{H^{k+1}\times H^k},
\qquad
\|\hat{\phi}_h-\phi\|_1 \le \bar{\alpha}_{31} h^{k} \|\phi\|_{H^{k +1}}.
\end{align}
(ii)~Suppose Hypothesis~\ref{hyp:L2} additionally holds.
Then, there exist positive constants~$\alpha_{32}$ and~$\bar{\alpha}_{32}$ independent of $h$ such that for any~$h$
\begin{align}
\|\hat{w}_h - w\|_0 \le \alpha_{32} h^{k+1} \|(w,r)\|_{H^{k+1}\times H^k},
\qquad
\|\hat{\phi}_h - \phi\|_0 \le \bar{\alpha}_{32} h^{k+1} \|\phi\|_{H^{k+1}}.
\label{ieq:StokesPoissonPrj_L2}
\end{align}
\end{subequations}
\end{Prop}
\subsection{An estimate at each time step}
Let $(\hat{u}_h, \hat{p}_h, \hat{\theta}_h)(t) \equiv \Pi_h^{\rm SP} (u, p, \theta)(t) \in V_h\times Q_h\times \Psi_h$ for $t\in [0,T]$, and let
\begin{align*}
e_{uh}^n \equiv u_h^n-\hat{u}_h^n,
\quad
e_{ph}^n \equiv p_h^n-\hat{p}_h^n, 
\quad
e_{\theta h}^n \equiv \theta_h^n-\hat{\theta}_h^n,
\quad
\eta_u (t) \equiv (u-\hat{u}_h)(t), 
\quad
\eta_\theta (t) \equiv (\theta-\hat{\theta}_h)(t).
\end{align*}
Then, from~\eqref{scheme_NCP}, \eqref{eq:StokesPoissonPrj} and~\eqref{eq:NCP_weak}, we have for $n\ge 1$
\begin{align}
(\ol{D}_{\Delta t}e_{uh}^n, v_h) + (\ol{D}_{\Delta t}e_{\theta h}^n, \psi_h) + \mathcal{A}_h \bigl( (e_{uh}^n,e_{ph}^n), (v_h,q_h) \bigr) + a_\theta (e_{\theta h}^n, \psi_h) = \lA R_{uh}^n, v_h\rA + \lA R_{\theta h}^n, \psi_h\rA, \notag\\
\forall (v_h, q_h, \psi_h)\in V_h\times Q_h\times \Psi_h,
\label{eq:error}
\end{align}
where
\begin{align*}
R_{uh}^n & \equiv \sum_{i=1}^7 R_{uhi}^n, & R_{\theta h}^n & \equiv \sum_{i=1}^4 R_{\theta hi}^n, \\
R_{uh1}^n & \equiv \fz{Du^n}{Dt} - \fz{u^n - u^{n-1} \circ X_1(u^{n-1},\Delta t)}{\Delta t}, &
R_{\theta h1}^n & \equiv \fz{D\theta^n}{Dt} - \fz{\theta^n - \theta^{n-1} \circ X_1(u^{n-1},\Delta t)}{\Delta t}, \\
R_{uh2}^n & \equiv \fz{1}{\Delta t}\Bigl\{ u^{n-1}\circ X_1(u_h^{n-1},\Delta t) - u^{n-1} \circ X_1(u^{n-1},\Delta t) \Bigr\}, &
R_{\theta h2}^n & \equiv \fz{1}{\Delta t}\Bigl\{ \theta^{n-1}\circ X_1(u_h^{n-1},\Delta t) - \theta^{n-1} \circ X_1(u^{n-1},\Delta t) \Bigr\}, \\
R_{uh3}^n & \equiv \fz{1}{\Delta t}\Bigl\{ \eta_u^n - \eta_u^{n-1} \circ X_1(u_h^{n-1},\Delta t) \Bigr\}, &
R_{\theta h3}^n & \equiv \fz{1}{\Delta t}\Bigl\{ \eta_\theta^n - \eta_\theta^{n-1} \circ X_1(u_h^{n-1},\Delta t) \Bigr\},\\
R_{uh4}^n & \equiv -\fz{1}{\Delta t} \Bigl\{ e_{uh}^{n-1} - e_{uh}^{n-1} \circ X_1(u_h^{n-1},\Delta t) \Bigr\}, &
R_{\theta h4}^n & \equiv -\fz{1}{\Delta t} \Bigl\{ e_{\theta h}^{n-1} - e_{\theta h}^{n-1} \circ X_1(u_h^{n-1},\Delta t) \Bigr\},\\
R_{uh5}^n & \equiv - \eta_\theta^{n-1} \beta^n, 
\quad
R_{uh6}^n \equiv e_{\theta h}^{n-1} \beta^n,
\quad
R_{uh7}^n \equiv -(\theta^n - \theta^{n-1}) \beta^n.
\end{align*}
We note that
\begin{align}
(e_{uh}^0, e_{ph}^0, e_{\theta h}^0) = (u_h^0, p_h^0, \theta_h^0) - (\hat{u}_h^0, \hat{p}_h^0, \hat{\theta}_h^0) = \Pi_h^{\rm SP} (u^0, 0, \theta^0) - \Pi_h^{\rm SP} (u^0, p^0, \theta^0).
\label{eq:initial_approx_value}
\end{align}
\begin{Prop}\label{prop:eh_epsh_Gronwall}
(i)~Let $(u, p, \theta)^0 \in (V\times Q\times \Psi) \cap \mathbb{H}^{k+1}$ be given.
Assume $\nabla\cdot u^0 =0$.
Then, there exists a positive constant $c_I$ independent of $h$ such that for any $h$
\begin{align}
\sqrt{\nu} \|D(e_{uh}^0)\|_0+\sqrt{\fz{\delta_0}{2}} |e_{ph}^0|_h+\sqrt{\fz{\kappa}{2}} \|\nabla e_{\theta h}^0\|_0 \le c_I h^k.
\label{ieq:eh0_epsh0}
\end{align}
(ii)~Let $n\in \{1,\cdots,N_T\}$ be a fixed number and let $(u_h^{n-1}, \theta_h^{n-1}) \in V_h \times \Psi_h$ be known.
Suppose the inequality
\begin{align}
\Delta t |u_h^{n-1}|_{1,\infty} \le 1/4
\label{ieq:dt_uh_1inf}
\end{align}
holds.
Then, there exists a unique solution~$(u_h^n,p_h^n,\theta_h^n)\in V_h \times Q_h \times \Psi_h$ of~\eqref{scheme_NCP}.
\medskip\\
(iii)~Furthermore, suppose Hypothesis~\ref{hyp:regularity} and the inequality
\begin{align}
\Delta t |u|_{C(W^{1,\infty})} \le 1/4
\label{ieq:dt_u_1inf}
\end{align}
hold.
Let $p_h^{n-1}\in Q_h$ be known and suppose the equation
\begin{align}
b(u_h^{n-1},q_h) - \delta_0 \mathcal{C}_h(p_h^{n-1},q_h) = 0,\quad \forall q_h\in Q_h,\label{eq:b_Ch_n-1}
\end{align}
holds.
Then, it holds that
\begin{align}
&\ol{D}_{\Delta t} \biggl( \nu\|D(e_{uh}^n)\|_0^2 +\fz{\delta_0}{2}|e_{ph}^n|_h^2 +\fz{\kappa}{2}\|\nabla e_{\theta h}^n\|_0^2 \biggr) + \fz{1}{2} (\|\ol{D}_{\Delta t}e_{uh}^n\|_0^2 + \|\ol{D}_{\Delta t}e_{\theta h}^n\|_0^2)
\le A_1(\|u_h^{n-1}\|_{0,\infty}) \biggl( \nu\|D(e_{uh}^{n-1})\|_0^2 +\fz{\kappa}{2}\|\nabla e_{\theta h}^{n-1}\|_0^2 \biggr) \notag\\
&\qquad
+ A_2(\|u_h^{n-1}\|_{0,\infty}) \biggl\{ \Delta t \biggl( \|u\|_{Z^2(t^{n-1},t^n)}^2 + \|\theta\|_{Z^2(t^{n-1},t^n)}^2 \biggr)
+ h^{2k} \Bigl( \fz{1}{\Delta t} \|(u,p,\theta)\|_{H^1(t^{n-1},t^n; \mathbb{H}^{k+1})}^2 +1 \Bigr) \Bigr\},
\label{ieq:eh_epsh_Gronwall}
\end{align}
where $A_i$, $i=1,2$, are functions defined by
\begin{align*}
A_i(\xi)\equiv c_i (\xi^2 +1 )
\end{align*}
and $c_i$, $i=1, 2$, are positive constants independent of $h$ and $\Delta t$.
They are defined by~\eqref{def:c1_c2} below.
\end{Prop}
\par
For the proof we use the next lemma, whose proof is omitted here.
We note that the proofs of estimates of~$\|R_{uhi}\|_0,\ i=1, \cdots, 4,$ in the case $k = 1$ are given in~\cite{NT-2015-M2AN}, 
that those in the case $k = 2$ are similarly obtained by Proposition~\ref{prop:StokesPoissonPrj},
that the proofs of the estimates of~$\|R_{\theta hi}\|_0,\ i=1, \cdots, 4,$ are similar 
and that the proofs of the estimates of~$\|R_{uhi}\|_0,\ i=5,6,7,$ are easier than them.
\begin{Lem}\label{lem:estimates_R}
Suppose Hypothesis~\ref{hyp:regularity} holds.
Let $n\in \{1,\cdots,N_T\}$ be a fixed number and let $u_h^{n-1} \in V_h$ be known.
Then, under the conditions \eqref{ieq:dt_uh_1inf} and \eqref{ieq:dt_u_1inf} it holds that
\begin{subequations}
\begin{align}
\| R_{uh1}^n \|_0 & \le c_u\sqrt{\Delta t} \|u\|_{Z^2(t^{n-1},t^n)}, 
&
\| R_{\theta h1}^n \|_0 & \le c_{(u,\theta)} \sqrt{\Delta t} \|\theta\|_{Z^2(t^{n-1},t^n)},
\label{ieq:R1}\\
\| R_{uh2}^n \|_0 & \le c_u \bigl( \|e_{uh}^{n-1}\|_0 + h^k \|(u,p)^{n-1}\|_{H^{k+1}\times H^k} \bigr), 
&
\| R_{\theta h2}^n \|_0 & \le c_\theta \bigl( \|e_{uh}^{n-1}\|_0 + h^k \|(u,p)^{n-1}\|_{H^{k+1}\times H^k} \bigr), 
\label{ieq:R2}\\
\| R_{uh3}^n \|_0 & \le \fz{ch^k}{\sqrt{\Delta t}} (\|u_h^{n-1}\|_{0,\infty} + 1) \| (u, p) \|_{H^1(t^{n-1},t^n; H^{k+1}\times H^k)}, 
&
\| R_{\theta h3}^n \|_0 & \le \fz{ch^k}{\sqrt{\Delta t}} (\|u_h^{n-1}\|_{0,\infty} + 1) \| \theta \|_{H^1(t^{n-1},t^n; H^{k+1})}, 
\label{ieq:R3}\\
\| R_{uh4}^n \|_0 & \le c \|u_h^{n-1}\|_{0,\infty} \|e_{uh}^{n-1}\|_1,
&
\| R_{\theta h4}^n \|_0 & \le c \|u_h^{n-1}\|_{0,\infty} \|e_{\theta h}^{n-1}\|_1,
\label{ieq:R4}\\
\| R_{uh5}^n \|_0 & \le c h^k \|\theta^{n-1}\|_{k+1},
\qquad
\| R_{uh6}^n \|_0 \le c \|e_{\theta h}^{n-1}\|_0,
&
\| R_{uh7}^n \|_0 & \le c \sqrt{\Delta t} \|\theta\|_{H^1(t^{n-1},t^n; L^2)},
\label{ieq:R7}
\end{align}
\end{subequations}
\end{Lem}
{\it Proof of Proposition~\ref{prop:eh_epsh_Gronwall}.}
\ \ 
We prove~(i).
From~\eqref{eq:initial_approx_value} we have
\begin{align*}
\|D(e_{uh}^0)\|_0 & \le \|e_{uh}^0\|_1 = \|u_h^0 - \hat{u}_h^0\|_1 \le \|u_h^0 - u^0\|_1 +\| u^0 - \hat{u}_h^0\|_1 \le 2\alpha_{31} h^k \|(u, p)^0\|_{H^{k+1} \times H^k}, \\
\|\nabla e_{\theta h}^0\|_0 & = \|\nabla (\theta_h^0 - \hat{\theta}_h^0)\|_0 \le \|\theta_h^0 - \theta^0\|_1 + \| \theta^0 - \hat{\theta}_h^0\|_1 = 2 \| \hat{\theta}_h^0 - \theta^0 \|_1 \le 2 \bar{\alpha}_{31} h^k \|\theta^0\|_{k+1},
\end{align*}
and in the case~$\delta_0=1~(k=1)$,
\begin{align*}
|e_{ph}^0|_h 
& = |p_h^0-\hat{p}_h^0|_h \le | p_h^0 - 0 |_h + |\hat{p}_h^0-p^0|_h + |p^0|_h 
\le \alpha_{20} \bigl( \| p_h^0 - 0 \|_0 + \|\hat{p}_h^0 - p^0\|_0 \bigr) + h\|p^0\|_1 \\
& \le (2\alpha_{20}\alpha_{31} + 1) h \| (u, p)^0 \|_{H^{k+1}\times H^k},
\end{align*}
which imply~\eqref{ieq:eh0_epsh0} for $c_I \equiv \{ 2\sqrt{\nu}\alpha_{31} + \sqrt{\delta_0/2} (2\alpha_{20}\alpha_{31}+1) \} \, \|(u,p)^0\|_{H^{k+1}\times H^k} + \sqrt{2\kappa}\bar{\alpha}_{31} \|\theta^0\|_{k+1}$.
\par
(ii) is obtained from~\eqref{ieq:dt_uh_1inf} and Remark~\ref{rmk:scheme}.
\par
We prove~(iii).
Substituting $(\ol{D}_{\Delta t}e_{uh}^n, 0, \ol{D}_{\Delta t}e_{\theta h}^n)$ into $(v_h, q_h, \psi_h)$ in~\eqref{eq:error}, we have
\begin{align}
& \|\ol{D}_{\Delta t} e_{uh}^n\|_0^2 + \|\ol{D}_{\Delta t} e_{\theta h}^n\|_0^2 + \ol{D}_{\Delta t} \biggl( \nu\|D(e_{uh}^n)\|_0^2 + \fz{\kappa}{2}\|\nabla e_{\theta h}^n\|_0^2 \biggr) + b(\ol{D}_{\Delta t} e_{uh}^n, e_{ph}^n) \le \| R_{uh}^n\|_0 \|\ol{D}_{\Delta t}e_{uh}^n\|_0 + \| R_{\theta h}^n\|_0 \|\ol{D}_{\Delta t}e_{\theta h}^n\|_0,
\label{proof_proposition_1}
\end{align}
where $X_1(u^{n-1},\Delta t)$ in $R_{uhi}^n$ and $R_{\theta hi}^n$, $i=1, 2$, maps $\Omega$ onto $\Omega$ by~\eqref{ieq:dt_u_1inf}.
From~\eqref{eq:b_Ch_n-1}, \eqref{scheme_NCP} with $(v_h, q_h, \psi_h) = (0, q_h, 0)\in V_h\times Q_h\times \Psi_h$ and $(\hat{u}_h^n, \hat{p}_h^n, \hat{\theta}_h^n)=\Pi_h^{\rm SP}(u^n, p^n, \theta^n)$ by~\eqref{eq:StokesPoissonPrj} we have
\begin{subequations}\label{eqns:uptheta_h_hat_uptheta_h_n-1_n}
\begin{align}
b(u_h^i, q_h) - \delta_0\mathcal{C}_h(p_h^i, q_h) &= 0, && \forall q_h \in Q_h,\\
b(\hat{u}_h^i, q_h) - \delta_0\mathcal{C}_h(\hat{p}_h^i, q_h) & = b(u^i, q_h) =0, && \forall q_h\in Q_h,
\end{align}
\end{subequations}
for $i=n-1$ and $n$.
The equalities~\eqref{eqns:uptheta_h_hat_uptheta_h_n-1_n} imply that
\begin{align*}
b(\ol{D}_{\Delta t} e_{uh}^n, q_h) - \delta_0\mathcal{C}_h(\ol{D}_{\Delta t} e_{ph}^n, q_h) = 0,\quad \forall q_h\in Q_h,
\end{align*}
which leads to
\begin{align}
-b(\ol{D}_{\Delta t} e_{uh}^n, e_{ph}^n) + \delta_0\mathcal{C}_h(\ol{D}_{\Delta t} e_{ph}^n, e_{ph}^n) = 0
\label{proof_proposition_2}
\end{align}
by putting $q_h=-e_{ph}^n\in Q_h$.
Adding~\eqref{proof_proposition_2} to~\eqref{proof_proposition_1} and using Lemma~\ref{lem:estimates_R},
we have
\begin{align*}
&\|\ol{D}_{\Delta t}e_{uh}^n\|_0^2 + \|\ol{D}_{\Delta t}e_{\theta h}^n\|_0^2 + \ol{D}_{\Delta t}\biggl( \nu\|D(e_{uh}^n)\|_0^2 +\fz{\delta_0}{2}|e_{ph}^n|_h^2 + \fz{\kappa}{2}\|\nabla e_{\theta h}^n\|_0^2 \biggr) \\
&\quad \le \fz{1}{2}( \|\ol{D}_{\Delta t}e_{uh}^n\|_0^2 + \|\ol{D}_{\Delta t}e_{\theta h}^n\|_0^2 ) + \fz{1}{2} \Bigl( \sum_{i=1}^7\| R_{uhi}^n\|_0^2 + \sum_{i=1}^4\| R_{\theta hi}^n \|_0^2 \Bigr) \\
&\quad \le \fz{1}{2}( \|\ol{D}_{\Delta t}e_{uh}^n\|_0^2 + \|\ol{D}_{\Delta t}e_{\theta h}^n\|_0^2 )
+ c_{(u,\theta)} \Bigl[ ( \|u_h^{n-1}\|_{0,\infty}^2 + 1 ) ( \|D(e_{uh}^{n-1})\|_0^2 + \|\nabla e_{\theta h}^{n-1}\|_0^2 ) + \Delta t ( \|u\|_{Z^2(t^{n-1},t^n)}^2 + \|\theta\|_{Z^2(t^{n-1},t^n)}^2 ) \\
&\qquad + h^{2k} \|(u, p, \theta)\|_{C(\mathbb{H}^{k+1})}^2 + (\|u_h^{n-1}\|_{0,\infty}^2+1)\fz{h^{2k}}{\Delta t} \|(u,p,\theta)\|_{H^1(t^{n-1},t^n; \mathbb{H}^{k+1})}^2 \Bigr],
\end{align*}
which implies that
\begin{align*}
& \ol{D}_{\Delta t}\biggl( \nu\|D(e_{uh}^n)\|_0^2 +\fz{\delta_0}{2}|e_{ph}^n|_h^2 + \fz{\kappa}{2}\|\nabla e_{\theta h}^n\|_0^2 \biggr) + \fz{1}{2} (\|\ol{D}_{\Delta t} e_{uh}^n\|_0^2 + \|\ol{D}_{\Delta t} e_{\theta h}^n\|_0^2) \\
& \quad \le \fz{c_{(u,\theta)}}{\min\{\nu,\kappa\}} (\|u_h^{n-1}\|_{0,\infty}^2 + 1) \Bigl(\nu\|D(e_{uh}^{n-1})\|_0^2 + \fz{\kappa}{2}\|\nabla e_{\theta h}^{n-1}\|_0^2\Bigr) \\
&\qquad + c_{(u,p,\theta)} \Bigl[ \Delta t (\|u\|_{Z^2(t^{n-1},t^n)}^2+\|\theta\|_{Z^2(t^{n-1},t^n)}^2) + ( \|u_h^{n-1}\|_{0,\infty}^2+1 ) h^{2k} \Bigl( \fz{1}{\Delta t} \|(u, p, \theta)\|_{H^1(t^{n-1},t^n; \mathbb{H}^{k+1})}^2 + 1 \Bigr) \Bigr].
\end{align*}
Putting
\begin{align}
c_1 \equiv \fz{c_{(u,\theta)}}{\min\{\nu,\kappa\}},
\qquad
c_2 \equiv c_{(u,p,\theta)},
\label{def:c1_c2}
\end{align}
we obtain~\eqref{ieq:eh_epsh_Gronwall}.
\qed
%
%
%
%
%
%
%
%
%
\subsection{Proof of Theorem~\ref{thm:main_results}}
The proof is performed by induction through three steps.\medskip\\
\textit{Step~1} (Setting $c_0$ and $h_0$):\ 
Let $c_I$ and $A_i$, $i=1,2$, be the constant and the functions in Proposition~\ref{prop:eh_epsh_Gronwall}, respectively.
Let $a_1$, $a_2$, $c_\ast$, $\bar{c}_\ast$ and $c_3$ be constants defined by
\begin{align*}
a_1 & \equiv A_1(\|u\|_{C(L^\infty)}+1),\quad a_2 \equiv A_2(\|u\|_{C(L^\infty)}+1), \\
c_\ast & \equiv \fz{\alpha_1}{\sqrt{\nu}} \, c_3,\quad \bar{c}_\ast \equiv \bar{\alpha}_1 \sqrt{\fz{2}{\kappa}} \, c_3,
\quad c_3 \equiv \exp(a_1T/2) \max \Bigl\{ a_2^{1/2} (\|u\|_{Z^2}+\|\theta\|_{Z^2}),  
a_2^{1/2} \bigl( \|(u,p,\theta)\|_{H^1(\mathbb{H}^{k+1})} + T^{1/2} \bigr) + c_I\Bigr\}.
\end{align*}
We can choose sufficiently small positive constants $c_0$ and $h_0$ such that
\begin{subequations}\label{def:c0_h0}
\begin{align}
\alpha_{21} \bigl\{ c_\ast (c_0h_0^{d/12}+h_0^{k-d/6}) + (\alpha_{24} + \alpha_{31}) h_0^{k-d/6}\|(u,p)\|_{C(H^{k+1}\times H^k)} \bigr\} & \le 1, 
\label{def:c0_h0_Linf}\\
c_0 \bigl[ \alpha_{22} \bigl\{ c_\ast (c_0+h_0^{k-d/4}) + (\alpha_{24} + \alpha_{31}) h_0^{k-d/4}\|(u,p)\|_{C(H^{k+1}\times H^k)} \bigr\} 
+ \alpha_{23} h_0^{d/4}\|u\|_{C(W^{1,\infty})} \bigr] & \le 1/4,
\label{def:c0_h0_W1inf}
\end{align}
\end{subequations}
since all the powers of $h_0$ are positive.
\medskip\\
\textit{Step~2} (Induction):\ 
For $n\in\{0,\cdots,N_T\}$ we define property P($n$) as follows:
\begin{align*}
\mbox{P($n$)}:
\left\{
\begin{aligned}
&
\!
\begin{aligned}
&
{\rm (a)}~\nu\|D(e_{uh}^n)\|_0^2 +\fz{\delta_0}{2}|e_{ph}^n|_h^2 +\fz{\kappa}{2}\|\nabla e_{\theta h}^n\|_0^2 + \fz{1}{2}(\|\ol{D}_{\Delta t}e_{uh}\|_{\ell^2_n(L^2)}^2 + \|\ol{D}_{\Delta t}e_{\theta h}\|_{\ell^2_n(L^2)}^2)
\\ & 
\qquad\quad
\le \exp (a_1 n\Delta t)
\Bigl[ \nu\|D(e_{uh}^0)\|_0^2 +\fz{\delta_0}{2}|e_{ph}^0|_h^2 + \fz{\kappa}{2}\|\nabla e_{\theta h}^0\|_0^2 \\
&  \qquad\qquad  
+ a_2\Bigl\{ \Delta t^2 (\|u\|_{Z^2(0,t^n)}^2 + \|\theta\|_{Z^2(0,t^n)}^2)
+ h^{2k} \bigl( \|(u,p,\theta)\|_{H^1(0,t^n; \mathbb{H}^{k+1})}^2 + n\Delta t \bigr)
\Bigr\} \Bigr], 
\end{aligned}
\\
& {\rm (b)}~\|u_h^n\|_{0,\infty} \le \|u\|_{C(L^\infty)}+1,\\
& {\rm (c)}~\ \Delta t |u_h^n|_{1,\infty} \le 1/4,
\end{aligned}
\right.
\end{align*}
where $\|\ol{D}_{\Delta t}e_{uh}\|_{\ell^2_n(L^2)} = \|\ol{D}_{\Delta t}e_{\theta h}\|_{\ell^2_n(L^2)} = 0$ for $n=0$.
P($n$)-(a) can be rewritten as
\begin{align}
x_n + \Delta t\sum_{i=1}^n y_i \le \exp( a_1n\Delta t ) \Bigl(x_0 + \Delta t \sum_{i=1}^n b_i \Bigr),
\label{ieq:proof_thm_P_n}
\end{align}
where
\begin{align*}
x_n & \equiv \nu\|D(e_{uh}^n)\|_0^2 +\fz{\delta_0}{2}|e_{ph}^n|_h^2 + \fz{\kappa}{2}\|\nabla e_{\theta h}^n\|_0^2,
\qquad
y_i \equiv \fz{1}{2} (\|\ol{D}_{\Delta t}e_{uh}^i\|_0^2+\|\ol{D}_{\Delta t}e_{\theta h}^i\|_0^2),\\
b_i & \equiv a_2\Bigl\{ \Delta t (\|u\|_{Z^2(t^{i-1},t^i)}^2 + \|\theta\|_{Z^2(t^{i-1},t^i)}^2) + h^{2k} \Bigl( \fz{1}{\Delta t} \|(u,p,\theta)\|_{H^1(t^{i-1},t^i; \mathbb{H}^{k+1})}^2 +1 \Bigr) \Bigr\}.
\end{align*}
\par
We firstly prove the general step in the induction.
Supposing that P($n-1$) holds true for an integer $n\in\{1,\cdots,N_T\}$, we prove that P($n$) also holds.
Since P($n-1$)-(c) is nothing but~\eqref{ieq:dt_uh_1inf}, there exists a unique solution $(u_h^n, p_h^n, \theta_h^n)\in V_h\times Q_h\times \Psi_h$ of equation~\eqref{scheme_NCP} from Proposition~\ref{prop:eh_epsh_Gronwall}-(ii).
We prove P($n$)-(a).
\eqref{ieq:dt_u_1inf} holds thanks to the estimate,
\begin{align*}
\Delta t |u|_{C(W^{1,\infty})} \le c_0h_0^{d/4} |u|_{C(W^{1,\infty})} \le c_0 \alpha_{23} h_0^{d/4} |u|_{C(W^{1,\infty})} \le 1/4,
\end{align*}
from condition~\eqref{condition:h_dt}, Remark~\ref{rmk:interpolation_property} and~\eqref{def:c0_h0_W1inf}.
\eqref{eq:b_Ch_n-1} is obtained from~\eqref{scheme_NCP} for $n\ge 2$ and from $(u_h^0, p_h^0, \theta_h^0) = \Pi_h^{\rm SP}(u^0, 0, \theta^0)$ for $n=1$.
Hence \eqref{ieq:eh_epsh_Gronwall} holds from Proposition~\ref{prop:eh_epsh_Gronwall}-(iii).
Since the inequalities $A_i(\|u_h^{n-1}\|_{0,\infty}) \le a_i~(i=1,2)$ hold from P($n-1$)-(b),
\eqref{ieq:eh_epsh_Gronwall} implies
\begin{align*}
& \ol{D}_{\Delta t} x_n + y_n \le a_1 x_{n-1} + b_n, 
\end{align*}
which leads to
\begin{align}
x_n + \Delta t y_n \le \exp (a_1 \Delta t ) (x_{n-1} + \Delta t b_n)
\label{proof_thm1_1}
\end{align}
by $1 \le 1+a_1\Delta t \le \exp (a_1 \Delta t)$.
From~\eqref{proof_thm1_1} and P($n-1$)-(a)
we have that
\begin{align*}
x_n + \Delta t\sum_{i=1}^n y_i 
& \le \exp (a_1\Delta t) (x_{n-1} + \Delta t b_n) + \Delta t\sum_{i=1}^{n-1} y_i 
\le \exp (a_1\Delta t) \Bigl( x_{n-1} + \Delta t\sum_{i=1}^{n-1} y_i + \Delta t b_n \Bigr) \\
& \le \exp (a_1\Delta t) \biggl[ \exp\bigl\{ a_1(n-1)\Delta t\bigr\} \biggl(x_0 + \Delta t \sum_{i=1}^{n-1}b_i \biggr) + \Delta t b_n \biggr]
\le \exp (a_1n\Delta t) \biggl(x_0 + \Delta t \sum_{i=1}^n b_i \biggr).
\end{align*}
Thus, we obtain~P($n$)-(a).
\par
For the proofs of P($n$)-(b) and (c) we prepare the estimate of $\|e_{uh}^n\|_1$.
From P($n$)-(a) and~\eqref{ieq:eh0_epsh0} we have that
\begin{align}
& \nu\|D(e_{uh}^n)\|_0^2 +\fz{\delta_0}{2}|e_{ph}^n|_h^2 +\fz{\kappa}{2}\|\nabla e_{\theta h}^n\|_0^2 + \fz{1}{2} (\|\ol{D}_{\Delta t}e_{uh}\|_{\ell^2_n(L^2)}^2 + \|\ol{D}_{\Delta t}e_{\theta h}\|_{\ell^2_n(L^2)}^2) \notag\\
& \quad \le \exp ( a_1T ) \Bigl[ c_I^2 h^{2k} 
+ a_2\Bigl\{ \Delta t^2(\|u\|_{Z^2}^2+\|\theta\|_{Z^2}^2) + h^{2k} \bigl( \|(u,p,\theta)\|_{H^1(\mathbb{H}^{k+1})}^2 + T \bigr) \Bigr\} \Bigr] \notag\\
& \quad \le \exp(a_1T) \Bigl[ a_2\Delta t^2(\|u\|_{Z^2}^2+\|\theta\|_{Z^2}^2) + h^{2k} \Bigl\{ a_2 \bigl( \|(u,p,\theta)\|_{H^1(\mathbb{H}^{k+1})}^2 + T \bigr)+c_I^2\Bigr\} \Bigr] \notag\\
& \quad \le \bigl\{ c_3 (\Delta t +h^k) \bigr\}^2.
\label{ieq:c3}
\end{align}
\eqref{ieq:c3} implies
\begin{align}
\|e_{uh}^n\|_1 & \le \alpha_1 \|D(e_{uh}^n)\|_0 \le \fz{\alpha_1}{\sqrt{\nu}} c_3 (\Delta t + h^k) = c_\ast (\Delta t + h^k).
\label{ieq:ehn_H1}
\end{align}
\par
We prove P($n$)-(b) and~(c) as follows:
\begin{align*}
\|u_h^n\|_{0,\infty} &\le \|u_h^n-\Pi_hu^n\|_{0,\infty} + \|\Pi_hu^n\|_{0,\infty} \le \alpha_{21} h^{-d/6}\|u_h^n-\Pi_hu^n\|_1 + \|\Pi_hu^n\|_{0,\infty} \\
& \le \alpha_{21} h^{-d/6} (\|u_h^n-\hat{u}_h^n\|_1 + \|\hat{u}_h^n-u^n\|_1 + \|u^n-\Pi_hu^n\|_1) + \|\Pi_hu^n\|_{0,\infty} \\
& \le \alpha_{21} h^{-d/6} \{ c_\ast (\Delta t+h^k) + \alpha_{31} h^k\|(u, p)^n\|_{H^{k+1}\times H^k} + \alpha_{24} h^k \|u^n\|_{k+1} \} + \|u^n\|_{0,\infty} \\
& \le \alpha_{21} \{ c_\ast (c_0h_0^{d/12}+h_0^{k-d/6}) + (\alpha_{24} + \alpha_{31}) h_0^{k-d/6} \|(u, p)\|_{C(H^{k+1}\times H^k)} \} + \|u\|_{C(L^\infty)} \\
& \le 1 + \|u\|_{C(L^\infty)},\\
\Delta t |u_h^n|_{1,\infty} & \le c_0h^{d/4} ( \|u_h^n-\Pi_hu^n\|_{1,\infty}+\|\Pi_hu^n\|_{1,\infty} )
\le c_0h^{d/4} ( \alpha_{22} h^{-d/2}\|u_h^n-\Pi_hu^n\|_1+\|\Pi_hu^n\|_{1,\infty} )\\
& \le c_0 \{ \alpha_{22} h^{-d/4} (\|u_h^n-\hat{u}_h^n\|_1 + \|\hat{u}_h^n-u^n\|_1 + \|u^n-\Pi_hu^n\|_1 ) +h^{d/4}\|\Pi_hu^n\|_{1,\infty} \}\\
& \le c_0 [ \alpha_{22} h^{-d/4} \{ c_\ast (\Delta t+h^k) + \alpha_{31}h^k \|(u, p)^n\|_{H^{k+1}\times H^k} + \alpha_{24} h^k \|u^n\|_{k+1} \} + \alpha_{23} h^{d/4}\|u^n\|_{1,\infty} ]\\
& \le c_0 [ \alpha_{22} h^{-d/4} \{ c_\ast (c_0h^{d/4}+h^k) + (\alpha_{24} + \alpha_{31}) h^k \|(u, p)^n\|_{H^{k+1}\times H^k} \} + \alpha_{23} h^{d/4}\|u^n\|_{1,\infty} ]\\
& \le c_0 [ \alpha_{22} \{ c_\ast (c_0+h_0^{k-d/4}) + (\alpha_{24} + \alpha_{31}) h_0^{k-d/4} \|(u, p)\|_{C(H^{k+1}\times H^k)} \} + \alpha_{23} h_0^{d/4}\|u\|_{C(W^{1,\infty})} ] \\
& \le 1/4,
\end{align*}
from~\eqref{ieq:ehn_H1}, \eqref{condition:h_dt} and~\eqref{def:c0_h0}.
Therefore, P($n$) holds true.
\par
The proof of P($0$) is easier than that of the general step.
P($0$)-(a) obviously holds with equality.
P($0$)-(b) and (c) are obtained as follows:
\begin{align*}
\|u_h^0\|_{0,\infty} & \le \|u_h^0-\Pi_hu^0\|_{0,\infty} + \|\Pi_hu^0\|_{0,\infty} 
\le \alpha_{21} h^{-d/6} ( \|u_h^0 - u^0 \|_1 + \|u^0 - \Pi_hu^0\|_1 ) + \|\Pi_hu^0\|_{0,\infty} \\
& \le \alpha_{21} (\alpha_{31} + \alpha_{24}) h^{k-d/6} \|(u, p)^0 \|_{H^{k+1}\times H^k} + \|u^0\|_{0,\infty} \\
& \le 1 + \|u\|_{C(L^\infty)},\\
\Delta t |u_h^0|_{1,\infty} & \le c_0h^{d/4} ( \|u_h^0-\Pi_hu^0\|_{1,\infty}+\|\Pi_hu^0\|_{1,\infty} )
\le c_0h^{d/4} ( \alpha_{22} h^{-d/2} \|u_h^0-\Pi_hu^0\|_1 +\|\Pi_hu^0\|_{1,\infty} ) \\
& \le c_0 \{  \alpha_{22} h^{-d/4} ( \|u_h^0 - u^0\|_1 + \|u^0-\Pi_hu^0\|_1 ) + h^{d/4}\|\Pi_hu^0\|_{1,\infty} \} \\
& \le c_0 \{  \alpha_{22} (\alpha_{31} + \alpha_{24}) h^{k-d/4} \|(u, p)^0\|_{H^{k+1}\times H^k} + \alpha_{23} h^{d/4}\|u^0\|_{1,\infty} \} \\
& \le 1/4.
\end{align*}
Thus, the induction is completed.
\medskip\\
\textit{Step~3}:\ 
Finally we derive the results~(i), (ii) and (iii) of the theorem.
Since P($N_T$) holds true, we have~(i), (ii) and the estimates
\begin{align}
\|e_{\theta h}\|_{\ell^\infty(H^1)} \le \bar{\alpha}_1 \|\nabla e_{\theta h}^n\|_{\ell^\infty(L^2)} \le  \bar{c}_\ast (\Delta t + h^k),
\qquad
\|\ol{D}_{\Delta t}e_{uh}\|_{\ell^2(L^2)},\ \|\ol{D}_{\Delta t}e_{\theta h}\|_{\ell^2(L^2)} \le \sqrt{2} c_3 (\Delta t + h^k),
\label{ieq:bar_eh_H1_Deh_L2}
\end{align}
from~\eqref{ieq:c3}.
The first and second inequalities of~\eqref{ieq:main_results} in~(iii) are obtained from~\eqref{ieq:ehn_H1}, \eqref{ieq:bar_eh_H1_Deh_L2} and the estimates
\begin{align*}
\|u_h-u\|_{\ell^\infty(H^1)} & \le \|e_{uh}\|_{\ell^\infty(H^1)} + \|\eta_u\|_{\ell^\infty(H^1)}
\le \|e_{uh}\|_{\ell^\infty(H^1)} + \alpha_{31} h^k \|(u, p)\|_{C(H^{k+1}\times H^k)}, \\
\Bigl\| \ol{D}_{\Delta t}u_h^n - \prz{u^n}{t} \Bigr\|_0 
& \le \|\ol{D}_{\Delta t}e_{uh}^n\|_0 + \|\ol{D}_{\Delta t}\eta_u^n\|_0 + \Bigl\| \ol{D}_{\Delta t}u^n - \prz{u^n}{t} \Bigr\|_0 \\
& \le \|\ol{D}_{\Delta t}e_{uh}^n\|_0 + \fz{\alpha_{31} h^k}{\sqrt{\Delta t}}\|(u,p)\|_{H^1(t^{n-1},t^n; H^{k+1}\times H^k)} + \sqrt{\fz{\Delta t}{3}} \Bigl\| \prz{^2u}{t^2} \Bigr\|_{L^2(t^{n-1},t^n; L^2)}.
\end{align*}
The proofs of the forth and fifth inequalities of~\eqref{ieq:main_results} are similar.
\par
We prove the third inequality of~\eqref{ieq:main_results}.
From~\eqref{ieq:generalized_inf_sup}, \eqref{eq:error} with $(v_h,q_h,\psi_h)=(v_h,q_h,0)$, Lemma~\ref{lem:estimates_R} and~\eqref{ieq:stability} we have
\begin{align*}
& \|e_{ph}^n\|_0 \le \|(e_{uh}^n,e_{ph}^n)\|_{V\times Q} \le \fz{1}{\gamma^\ast} \sup_{(v_h,q_h)\in V_h\times Q_h} \fz{\mathcal{A}_h\bigl( (e_{uh}^n,e_{ph}^n), (v_h,q_h) \bigr)}{\|(v_h,q_h)\|_{V\times Q}} 
= \fz{1}{\gamma^\ast} \sup_{(v_h,q_h)\in V_h\times Q_h} \fz{\lA R_{uh}^n, v_h\rA -(\ol{D}_{\Delta t}e_{uh}^n, v_h)}{\|(v_h,q_h)\|_{V\times Q}}
\notag\\
& 
\le c_{(u,p,\theta)} \Bigl\{ \sqrt{\Delta t} (\|u\|_{Z^2(t^{n-1},t^n)} + \|\theta\|_{H^1(t^{n-1},t^n; L^2)})
+ h^k\Bigl( \fz{1}{\sqrt{\Delta t}} \|(u,p,\theta)\|_{H^1(t^{n-1},t^n; \mathbb{H}^{k+1})} + 1 \Bigr) 
+ \|e_{uh}^{n-1}\|_1 + \|e_{\theta h}^{n-1}\|_0 + \|\ol{D}_{\Delta t}e_{uh}^n\|_0 \Big\}.
\end{align*}
Combining \eqref{ieq:ehn_H1} and~\eqref{ieq:bar_eh_H1_Deh_L2} with the above estimate, we obtain $\|e_{ph}\|_{\ell^2(L^2)} \le c_{(u,p,\theta)}(\Delta t+h^k)$, which leads to the result from
\begin{align*}
\phantom{MMMMMmm}
\|p_h-p\|_{\ell^2(L^2)} \le \|e_{ph}\|_{\ell^2(L^2)} + \|\hat{p}_h-p\|_{\ell^2(L^2)} 
\le \|e_{ph}\|_{\ell^2(L^2)} + \sqrt{T} \alpha_{31} h^k \|(u,p)\|_{C(H^{k+1}\times H^k)}.
\phantom{MMMMMm}
\qed
\end{align*}
%
%
%
%
%
%
%
\subsection{Proof of Theorem~\ref{thm:main_results_L2}}
For the proof we use the next lemma without proof, since the proofs of estimates in the lemma are similar to or easier than those of~\cite{NT-2015-M2AN}.
\begin{Lem}\label{lem:estimates_R_L2}
Suppose Hypotheses~\ref{hyp:regularity} and~\ref{hyp:L2} hold.
Let $n\in \{1,\cdots,N_T\}$ be a fixed number and $u_h^{n-1} \in V_h$ be known.
Then, under the conditions \eqref{ieq:dt_uh_1inf} and \eqref{ieq:dt_u_1inf} we have that 
\begin{subequations}\label{ieqs:R_L2}
\begin{align}
\| R_{uh2}^n \|_0 
& \le c_u \bigl( \|e_{uh}^{n-1}\|_0 + h^{k+1} \|(u,p)^{n-1}\|_{H^{k+1}\times H^k}\bigr),
\qquad
\| R_{\theta h2}^n \|_0 
\le c_\theta \bigl( \|e_{uh}^{n-1}\|_0 + h^{k+1} \|(u,p)^{n-1}\|_{H^{k+1}\times H^k}\bigr),
\label{ieq:R2_L2}\\
\| R_{uh3}^n \|_{V_h^\prime} 
& \le c_u \Bigl[ \|(u,p)^{n-1}\|_{H^{k+1}\times H^k} \Bigl\{ \|e_{uh}^{n-1}\|_0 + h^{k+1} (\|(u,p)^{n-1}\|_{H^{k+1}\times H^k} + 1) \Bigr\} + \fz{h^{k+1}}{\sqrt{\Delta t}} \|(u,p)\|_{H^1(t^{n-1},t^n; H^{k+1}\times H^k)} \Bigr],
\label{ieq:R3_L2}\\
\| R_{\theta h3}^n \|_{\Psi_h^\prime} 
& \le c_u \Bigl[ \|\theta^{n-1}\|_{k+1} \Bigl\{ \|e_{uh}^{n-1}\|_0 + h^{k+1} ( \|(u,p)^{n-1}\|_{H^{k+1}\times H^k} +1 ) \Bigr\} + \fz{h^{k+1}}{\sqrt{\Delta t}} \|\theta\|_{H^1(t^{n-1},t^n; H^{k+1})} \Bigr],
\label{ieq:bar_R3_L2}\\
\| R_{uh4}^n \|_{V_h^\prime} 
& \le c_u \bigl( 1+h^{-d/6}\|e_{uh}^{n-1}\|_1 \bigr) \bigl( \|e_{uh}^{n-1}\|_0 + h^{k+1} \|(u,p)^{n-1}\|_{H^{k+1}\times H^k} \bigr),
\label{ieq:R4_L2}\\
\| R_{\theta h4}^n \|_{\Psi_h^\prime} 
& \le c_u \bigl\{ \|e_{\theta h}^{n-1}\|_0 + h^{-d/6}\|e_{\theta h}^{n-1}\|_1 \bigl( \|e_{uh}^{n-1}\|_0 + h^{k+1} \|(u,p)^{n-1}\|_{H^{k+1}\times H^k} \bigr) \bigr\},
\label{ieq:bar_R4_L2}\\
\| R_{uh5}^n \|_0 & \le c h^{k+1} \|\theta^{n-1}\|_{k+1}.
\label{ieq:R5_L2}
\end{align}
\end{subequations}
\end{Lem}
{\it Proof of Theorem~\ref{thm:main_results_L2}.}
\ \ 
Since we have $\|e_{uh}\|_{\ell^\infty(H^1)}\le c_\ast (\Delta t+h^k) \le c_\ast (c_0+h_0^{k-d/4})h^{d/4}$ and $\|e_{\theta h}\|_{\ell^\infty(H^1)}\le \bar{c}_\ast (\Delta t+h^k) \le \bar{c}_\ast (c_0+h_0^{k-d/4})h^{d/4}$ from \eqref{ieq:ehn_H1}, \eqref{ieq:bar_eh_H1_Deh_L2} and \eqref{condition:h_dt}, \eqref{ieq:R4_L2} and~\eqref{ieq:bar_R4_L2} imply
\begin{align}
\| R_{uh4}^n \|_{V_h^\prime} & \le c_u c_\ast \bigl( \|e_{uh}^{n-1}\|_0 + h^{k+1} \|(u,p)^{n-1}\|_{H^{k+1}\times H^k} \bigr),
&
\| R_{\theta h4}^n \|_{\Psi_h^\prime} 
& \le c_u \bigl\{ \|e_{\theta h}^{n-1}\|_0 + \bar{c}_\ast ( \|e_{uh}^{n-1}\|_0 + h^{k+1} \|(u,p)^{n-1}\|_{H^{k+1}\times H^k}) \bigr\}.
\label{ieq:R4_L2_1}
\end{align}
Substituting $(e_{uh}^n, -e_{ph}^n, e_{\theta h}^n)$ into $(v_h,q_h,\psi_h)$ in~\eqref{eq:error} and using Lemma~\ref{lem:Korn}, \eqref{ieq:R1}, \eqref{ieq:R7}, \eqref{ieqs:R_L2}, \eqref{ieq:R4_L2_1} and $ab \le a^2/(4\beta)+\beta b^2~(\beta>0)$, we have
\begin{align*}
& \ol{D}_{\Delta t} \Bigl( \fz{1}{2}\|e_{uh}^n\|_0^2 + \fz{1}{2}\|e_{\theta h}^n\|_0^2 \Bigr) + \fz{2\nu}{\alpha_1^2} \|e_{uh}^n\|_1^2 + \delta_0|e_{ph}^n|_h^2 + \fz{\kappa}{\bar{\alpha}_1^2} \|e_{\theta h}^n\|_1^2 
\le \lA R_{uh}^n, e_{uh}^n \rA + \lA R_{\theta h}^n, e_{\theta h}^n \rA \\
& \le \Bigl( \sum_{i=1,2,5,6,7}\|R_{uhi}^n\|_0 + \sum_{i=3, 4} \|R_{uhi}^n\|_{V_h^\prime} \Bigr) \|e_{uh}^n\|_1 + \Bigl( \sum_{i=1,2}\|R_{\theta hi}^n\|_0 + \sum_{i=3, 4} \|R_{\theta hi}^n\|_{\Psi_h^\prime} \Bigr) \|e_{\theta h}^n\|_1 \\
& \le \fz{\nu}{\alpha_1^2} \|e_{uh}^n\|_1^2 + \fz{\kappa}{2\bar{\alpha}_1^2} \|e_{\theta h}^n\|_1^2 + \fz{\alpha_1^2}{4\nu} \Bigl( \sum_{i=1,2,5,6,7}\|R_{uhi}^n\|_0 + \sum_{i=3, 4} \|R_{uhi}^n\|_{V_h^\prime} \Bigr)^2 +  \fz{\bar{\alpha}_1^2}{4\kappa} \Bigl( \sum_{i=1,2}\|R_{\theta hi}^n\|_0 + \sum_{i=3, 4} \|R_{\theta hi}^n\|_{\Psi_h^\prime} \Bigr)^2 \\
& \le \fz{\nu}{\alpha_1^2} \|e_{uh}^n\|_1^2 + \fz{\kappa}{2\bar{\alpha}_1^2} \|e_{\theta h}^n\|_1^2 
+ c_{(u,\theta)} \Bigl[ ( 1 + c_\ast^2 + \bar{c}_\ast^2 + \|(u,p,\theta)\|_{C(\mathbb{H}^{k+1})}^2 ) \|e_{uh}^{n-1}\|_0^2 + \|e_{\theta h}^{n-1}\|_0^2 
+ \Delta t ( \|u\|_{Z^2(t^{n-1},t^n)}^2 + \|\theta\|_{Z^2(t^{n-1},t^n)}^2 ) \\
& \quad + \fz{h^{2(k+1)}}{\Delta t}\|(u,p,\theta)\|_{H^1(t^{n-1},t^n; \mathbb{H}^{k+1})}^2 + h^{2(k+1)} \Bigl\{ ( 1 + c_\ast^2 + \bar{c}_\ast^2 + \| (u,p,\theta) \|_{C(\mathbb{H}^{k+1})}^2) \| (u,p,\theta) \|_{C(\mathbb{H}^{k+1})}^2 \Bigr\} \Bigr],
\end{align*}
where $\|e_{uh}^n\|_0 \le \|e_{uh}^n\|_1$ and $\|e_{\theta h}^n\|_0 \le \|e_{\theta h}^n\|_1$ have been employed.
Hence, we have that
\begin{align*}
& \ol{D}_{\Delta t} \Bigl( \fz{1}{2}\|e_{uh}^n\|_0^2 + \fz{1}{2}\|e_{\theta h}^n\|_0^2 \Bigr) + \fz{\nu}{\alpha_1^2} \|e_{uh}^n\|_1^2 + \fz{\kappa}{2\bar{\alpha}_1^2} \|e_{\theta h}^n\|_1^2 \\
&
\le c_{(u,p,\theta)} \Bigl\{ \|e_{uh}^{n-1}\|_0^2 + \|e_{\theta h}^{n-1}\|_0^2
+ \Delta t (\|u\|_{Z^2(t^{n-1},t^n)}^2 + \|\theta\|_{Z^2(t^{n-1},t^n)}^2) +  h^{2(k+1)} \Bigl( \fz{1}{\Delta t}\|(u,p,\theta)\|_{H^1(t^{n-1},t^n; \mathbb{H}^{k+1})}^2 + 1 \Bigr) \Bigr\}.
\end{align*}
From discrete Gronwall's inequality there exists a positive constant~$c_4$ independent of $h$ and $\Delta t$ such that
\begin{align}
\|e_{uh}\|_{\ell^\infty(L^2)},\ \|e_{\theta h}\|_{\ell^\infty(L^2)} & \le c_4 ( \|e_{uh}^0\|_0 + \|e_{\theta h}^0\|_0 + \Delta t + h^{k+1} ).
\label{ieq:eh_Linf_L2}
\end{align}
\eqref{ieq:eh_Linf_L2} and the estimates
\begin{align*}
\|e_{uh}^0\|_0 &\le \|u_h^0-u^0\|_0 + \|u^0-\hat{u}_h^0\|_0 \le 2\alpha_{32} h^{k+1} \|(u,p)^0\|_{H^{k+1}\times H^k},\\
\|u_h-u\|_{\ell^\infty(L^2)} & \le \|e_{uh}\|_{\ell^\infty(L^2)} + \|\eta_u\|_{\ell^\infty(L^2)} \le \|e_{uh}\|_{\ell^\infty(L^2)} + \alpha_{32}h^{k+1} \|(u,p)\|_{C(H^{k+1}\times H^k)},
\end{align*}
imply the first inequality of~\eqref{ieq:L2}.
The second is obtained similarly.
\qed
%
%
%
%
%
\section{Conclusions}\label{sec:conclusions}
We have proved optimal error estimates of stable and stabilized LG schemes for natural convection problems under a condition~$\Delta t = O(h^{d/4})$.
The stable and stabilized schemes employ P2/P1/P2~($k=2$) and P1/P1/P1~($k=1$) finite elements, respectively.
Both schemes have common advantages of the LG method, i.e., robustness for convection-dominated problems and symmetry of the coefficient matrix of the system of linear equations, where the solution of each scheme is obtained by solving alternatively the equations of flow and temperature.
The stabilized scheme has an additional advantage, a small number of degrees of freedom, which leads to efficient computation especially for three-dimensional problems.
The proofs of the optimal error estimates have been done by extending the arguments employed for the proofs of error estimates of stable and stabilized LG schemes for the Navier-Stokes equations in~\cite{NT-2015-M2AN,Suli-1988}.
We note that more general elements can be also dealt with, cf. Remark~\ref{rmk:generalization}.
It is not difficult to prove the optimal error estimates of corresponding two-step stable and stabilized LG schemes of second-order in time by combining the argument in this paper with those of~\cite{BMMR-1997,NT-2015-CREST}, and the second-order numerical convergence of a two-step stable LG scheme proposed in~\cite{BB-2011} is ensured mathematically.
%
%
%
\section*{Acknowledgements}
This work was supported by JSPS (the Japan Society for the Promotion of Science) under the Japanese-German Graduate Externship (Mathematical Fluid Dynamics) and Grant-in-Aid for Scientific Research~(S), No.~24224004.
The authors are indebted to JSPS also for Grant-in-Aid for Young Scientists~(B), No.~26800091 to the first author and for Grant-in-Aid for Scientific Research~(C), No.~25400212 to the second author.
%
%
%
%
%
%
%
%
%
%
%

%
%
%
%
\end{document}